\documentclass[10pt,a4paper]{amsart}
\usepackage{amssymb}
\usepackage[all]{xy}

 \theoremstyle{plain}

\newtheorem{teo}{Theorem}
\newtheorem{lem}[teo]{Lemma}
\newtheorem{pro}[teo]{Proposition}
\newtheorem{cor}[teo]{Corollary}

\newtheorem*{teoa}{Theorem A}

\newtheorem*{teob}{Theorem B}

\theoremstyle{definition}

\newtheorem{eje}[teo]{Example}

\newcommand{\F}{\mathbb{F}}
\newcommand{\N}{\mathbb{N}}

\newcommand{\eu}[1]{\mathfrak{#1}}

\DeclareMathOperator{\kernel}{ker}
\DeclareMathOperator{\image}{im}
\DeclareMathOperator{\rst}{res}
\DeclareMathOperator{\mor}{mor}
\DeclareMathOperator{\Hom}{Hom}

\DeclareMathOperator{\Aut}{Aut}

\DeclareMathOperator{\res}{res}
\DeclareMathOperator{\Spec}{Spec}

\DeclareMathOperator{\bp}{\mathbf{p}}
\DeclareMathOperator{\Cp}{\ca{C}_{\bp}}

\newcommand{\ca}[1]{\mathcal{#1}}

\newcommand{\argu}{\hbox to 7truept{\hrulefill}}

\title{Cohomology, fusion and a p-nilpotency criterion}
\author{Jon Gonz\'alez-S\'anchez}

\address{Jon Gonz\'alez-S\'anchez\\
Universidad de Cantabria\\
Departamento de Matem\'aticas, Estad\'{\i}stica y Computaci\'on \\
Facultad de Ciencias, Avda. de los Castros, s/n \\
E-39071 Santander, Spain}

\email{jon.gonzalez@unican.es}

\date{\today}

\begin{document}

\begin{abstract}
Let $G$ be a finite group, $p$ a fix prime  and $P$ a Sylow $p$-subgroup of $G$. In this short note we prove that 
if $p$ is odd, $G$ is $p$-nilpotent if and only if $P$ controls fusion of cyclic groups of order $p$. For the case $p=2$, 
we show that $G$ is $p$-nilpotent if and only if $P$ controls fusion of cyclic groups of order $2$ and $4$.
\end{abstract}

\maketitle

\section{Introduction}

Throughout the text let $p$ denote a fix prime. 
Let $G$ be a finite group and $P$ a Sylow $p$-subgroup of $G$. 
We denote by $H^\bullet (G,\F_p)$ the mod $p$ cohomology algebra. It is well known that the restriction map in cohomology
\begin{equation}
\label{eq:resisom}
\xymatrix{
H^\bullet (G,\F_p)\ar@{^{(}->}[r]&H^\bullet (P,\F_p)
}
\end{equation}
is injective (see \cite[Proposition 4.2.2]{ev:cgr}). Suppose that $G$ is $p$-nilpotent, i.e., $P$ has a normal complement $N$ in $G$. In this situation 
the composition 
\begin{equation}
\label{eq:composition}
\xymatrix{
P\ar[r] & G\ar[r] & G/N,
}
\end{equation}
is an isomorphism. Therefore the composition
\begin{equation}
\xymatrix{
H^\bullet (G/N,\F_p)\ar[0,2]^{inf^G_{G/N}} & & H^\bullet (G,\F_p)\ar[0,2]^{res^G_P} & & H^\bullet (P,\F_p)
}
\end{equation}
is also an isomorphism. This together with \eqref{eq:resisom} implies that, if 
$G$ is $p$-nilpotent, then the 
restriction map in cohomology
$\res^G_P :H^\bullet (G,\F_p)\longrightarrow H^\bullet (P,\F_p)$ is an isomorphism. 
The following result of M. Atiyah shows that the converse is also true.

\begin{teo}[Atiyah]
If $\res^G_P :H^i (G,\F_p)\longrightarrow H^i (P,\F_p)$ are isomorphisms for all $i$ big enough, 
then $G$ is $p$-nilpotent. In particular $G$ is $p$-nilpotent if and only if 
$\res^G_P :H^\bullet (G,\F_p)\longrightarrow H^\bullet (P,\F_p)$ is an isomorphism.
\end{teo}

\begin{proof}
A proof of this can be found in the introduction of \cite{qu:pnil}. 
\end{proof}

Atiyah $p$-nilpotency criterion uses the cohomology in high dimension.
Another cohomological criterion for $p$-nilpotency using cohomology in dimension $1$ 
was provided by J. Tate (\cite{tat:nil}).

\begin{teo}[Tate]
If $\res^G_P:H^1 (G,\F_p)\longrightarrow H^1 (P,\F_p)$ is an isomorphism, then $G$ is $p$-nilpotent.
\end{teo}

\begin{proof}
See  \cite{tat:nil}.
\end{proof}

D. Quillen generalized Atiyah's $p$-nilpotency criterion for odd primes (\cite{qu:pnil}). 

\begin{teo}[Quillen]
\label{teo:quipnil}
Let $p$ be an odd prime. Then $G$ is $p$-nilpotent if and only if $\res^G_P:H^\bullet (G,\F_p)\longrightarrow H^\bullet  (P,\F_p)$ is an $F$-isomorphism. 
\end{teo}

\begin{proof}
See \cite{qu:pnil}.
\end{proof}

Atiyah's $p$-nilpotency criterion can be reinterpreted in terms of $p$-fusion. 
We recall that  a subgroup $H$ of $G$ \textit{controls $p$-fusion in $G$} if 
\begin{itemize}
\item[(a)] $H$ contains a Sylow $p$-subgroup of $G$ and
\item[(b)] for any subgroup $A$ of $G$ and for any $g\in G$ such that $A, A^g\leq H$, 
there exists $x\in H$ such that for all $a\in A$, $a^g=a^x$.
\end{itemize}
By a result of G. Mislin \cite{mis:modp}, a subgroup $H$ of $G$ controls $p$-fusion in $G$ 
if and only if $\res^G_P:H^\bullet (G,\F_p)\longrightarrow H^\bullet  (H,\F_p)$ is an isomorphism. 
Using Mislin's result Atiyah's $p$-nilpotency criterion follows from Frobenious $p$-nilpotency criterion.

Mislin's type of result can also be provided 
for the concept of $F$-isomorphism.  In order to do this we introduce the following concept. Let $\ca{C}$ be a class of finite $p$-groups. We say that a subgroup $H$ of $G$ \textit{controls fusion of $\ca{C}$-groups in $G$} if
\begin{itemize}
\item[(a)] Any $\ca{C}$-subgroup of $G$ is conjugated to a subgroup of $H$ and
\item[(b)] for any $\ca{C}$-subgroup $A$ of $G$ and for any $g\in G$ such that $A, A^g\leq H$, there exists $x\in H$ such that for all $a\in A$, $a^g=a^x$.
\end{itemize}
The condition (b) can be rewritten as 
\begin{itemize}
\item[(b$^\prime$)] if $A$ is a $\ca{C}$-subgroup of $H$ and $g\in G$ satisfies that $A^g\leq H$, then $g\in C_G(A).H$. 
\end{itemize}

Theorem A bellow, which will be proved in Section \ref{s:cofus}, follows
naturally from Quillen's work on cohomology (see \cite{qu:strat} and \cite{qu:pnil}). Note that the 
``if'' was proved in \cite{GSW} and it is a direct consequence of Quillen's stratification (\cite{qu:strat}). 
The converse follows from a careful reading of \cite[Section 2]{qu:pnil}.

\begin{teoa}
Let $G$ be a finite group and $H$ a subgroup of $G$. Then $\res^G_H:H^\bullet (G,\F_p)\longrightarrow  H^\bullet (H,\F_p)$ 
is an $F$-isomorphims if and only if $H$ controls fusion of elementary abelian $p$-subgroups of $G$.
\end{teoa}

In Section \ref{s:pnilcriterion} we will prove the following $p$-nilpotency criterion that can be seen as a generalization of Quillen $p$-nilpotency criterion (Theorem \ref{teo:quipnil} above) to the prime $p=2$.

\begin{teob}
Let $G$ be a finite group and $P$ a Sylow $p$-subgroup of $G$. Then  the following two conditions are equivalent
\begin{enumerate}
\item $G$ is $p$-nilpotent.
\item $P$ controls fusion of cyclic subgroups of order $p$ in case $p$ is odd, and cyclic subgroups of order $2$ and $4$ in case $p=2$.
\end{enumerate}
\end{teob}

Note that Theorem A and Theorem B imply Quillen's $p$-nilpotency criterion. We will finish this short note by giving two applications of Theorem B. The first application will consist on reproving a result of H-W. Henn and S. Priddy  that implies that "most" finite groups are $p$-nilpotent (see \cite{hp:pcen}). The second application is a generalization to the prime $p=2$ of the following fact: if all elements of order $p$ of a finite group $G$ are in some upper center of $G$ and $p$ is an odd prime, then $G$ is $p$-nilpotent (see \cite{th:pcenp} and \cite{GSW}). For the prime $p=2$ we will show that if all elements of order $2$ and $4$ are in some upper center of $G$, then $G$ is $2$-nilpotent.

We would like to end this introduction with an example of Quillen \cite{qu:pnil} where 
the necessity of considering cyclic groups of order $2$ and $4$ for the case $p=2$ in Theorem B
is illustrated.

\begin{eje}
Consider $Q=\{ 1,-1,i,-i,j,-j,k,-k\}$ the quaternion group and $\alpha$ an automorphism of order $3$ that permutes $i$, $j$ and $k$. Let $G$ be the semidirect product between $Q$ and $\langle \alpha\rangle$ given by the action of $\alpha$ in $Q$. $A=\{ 1,-1\}$ is the only subgroup of exponent $2$ in $G$. Clearly $Q$ controls fusion of cyclic subgroup of order $2$. However $G$ is not $2$-nilpotent.
\end{eje}


\section{Cohomology and fusion}
\label{s:cofus}
The aim of this section is to sketch the proof of Theorem A. 
In subsections \ref{ss:strat}, \ref{ss:ea} and \ref{ss:spec} 
we will recall Quillen work in the mod $p$ cohomology algebra 
of a finite group. This will be used in subsection \ref{ss:cf} to 
prove Theorem A. 

For a finite group $G$ the {\it mod $p$ cohomology algebra}
\begin{equation}
\label{eq:bas-11}
H^\bullet(G)=H^\bullet(G,\F_p)
\end{equation}
is a finitely generated, connected, anti-commutative, $\N_0$-graded $\F_p$-algebra.

Let $\alpha_\bullet\colon A_\bullet\to B_\bullet$ be a homomorphism of
finitely generated, connected, anti-commutative, $\N_0$-graded $\F_p$-algebras.
Then $\alpha_\bullet$ is called an {\it $F$-isomorphism} if
$\kernel(\alpha_\bullet)$ is nilpotent, and for all $b\in B_n$ there exists
$k\geq 0$ such that $b^{p^k}\in\image(\alpha_\bullet)$.


\subsection{Quillen's stratification}
\label{ss:strat}
Let $G$ be a finite group. Let $\ca{E}_G$ denote the category
whose objects are the elementary abelian $p$-subgroups of $G$ and whose morphisms are given 
by conjugation, i.e., for $E, E^\prime\in \text{ob}(\ca{E}_G)$ one has
\begin{equation}
\label{eq:bas-4}
\mor_G(E,E^\prime)=\{\,i_g\colon E\to E^\prime\mid g\in G,\ g\,E\,g^{-1}\leq E^\prime\,\},
\end{equation}
where $i_g(e)= g\, e\, g^{-1}$, $e\in E$. Then
\begin{equation}
\label{eq:bas-5}
H^\bullet(\ca{E}_G)=\textstyle{\varprojlim_{\ca{E}_G} H^\bullet(E)}
\end{equation}
is a finitely generated, connected, anti-commutative, $\N_0$-graded $\F_p$-algebra.
Moreover, the restriction maps $\res^G_E$ yield a map
\begin{equation}
\label{eq:bas-6}
q_G=\textstyle{\prod_{E\in ob(\ca{E}_G)}\rst^G_E}\colon H^\bullet(G)\longrightarrow
H^\bullet(\ca{E}_G).
\end{equation}
The following result is  known as Quillen stratification. 

\begin{teo}[Quillen]
\label{thm:strat}
Let $G$ be a finite group. Then $q_G\colon H^\bullet(G)\to H^\bullet(\ca{E}_G)$
is an $F$-isomorphism.
\end{teo}

\begin{proof}
See \cite[Cor. 5.6.4]{ben:coh2} or \cite{qu:strat}.
\end{proof}


\subsection{Cohomology of elementary abelian $p$-groups} 
\label{ss:ea}
One can easily deduce the cohomology 
of an elementary abelian $p$-group from the cohomology of the cyclic group of exponent $p$ 
and the Kunneth formula.

\begin{lem}
Let $A$ be an elementary abelian $p$-group. Then
\begin{equation}
H^\bullet (A,\F_p)\cong \begin{cases}
\Lambda (A^*) \otimes S (\beta (A^* )) & \text{if $p$ is odd} \\
S (A^*) & \text{if $p=2$},
\end{cases}
\end{equation}
where $\Lambda$ denotes the exterior algebra functor, $S$ the symmetric algebra functor, 
$A^*=\Hom (A,\F_p)=H^1(A,\F_p)$ and $\beta$ the Bockstein homomorphism from $H^1(A,\F_p)$ to 
$H^2(A,\F_p)$.
\end{lem}

\begin{proof}
See \cite[Chap. 3 Section 5]{ev:cgr}.
\end{proof}

From the previous lemma one can easily deduces that

\begin{equation}
\label{eq:elemabelian}
H^\bullet (A,\F_p)/\sqrt{0}\cong 
S (A^*).
\end{equation}


\subsection{The spectrum of  $H(G)$}
\label{ss:spec}
Let $G$ be a finite group. Following  Quillen (\cite{qu:pnil}) we define
\begin{equation}
H(G)=\begin{cases}
\oplus_{i\geq 0} H^{2i}(G,\F_p) & \text{if $p$ is odd} \\
\oplus_{i\geq 0} H^{i}(G,\F_p)  & \text{if $p=2$}.
\end{cases}
\end{equation}
$H(G)$ is a graded commutative ring.  For an elementary abelian $p$-subgroup $A$ of $G$, 
denote by $\eu{g}_A$ 
the ideal of $H (G)$ consisting of elements $u$ such that $u|_A$ is nilpotent. From 
\eqref{eq:elemabelian}, $\res^G_A:H(G)\longrightarrow H (A)$ induces a 
monomorphism
\begin{equation}
\label{eq:spec-1}
\xymatrix{
H (G)/\eu{g}_A\ar@{^{(}->}[r]&S (A^*).
}
\end{equation}
In particular, the ideal  $\eu{g}_A$ is a prime ideal of $H(G)$. Furthermore,

\begin{teo}[Quillen]
\label{teo:conjugacion}
Let $A,A^\prime \subset G$ be elementary abelian subgroups of $G$. Then $\eu{g}_A\subseteq \eu{g}_{A^\prime}$ 
if and only if $A^\prime$ is conjugated to a subgroup of $A$. In particular $\eu{g}_A= \eu{g}_{A^\prime}$ 
if and only if $A$ and $A^\prime$ are conjugated in $G$.
\end{teo}

\begin{proof}
See \cite[Theorem 2.7]{qu:pnil}.
\end{proof}

Let us consider the extension of quotient fields associated to the monomorphism in \eqref{eq:spec-1},
\begin{equation}
\xymatrix{
k(\eu{g}_A)\ar@{^{(}->}[r]&k(A).
}
\end{equation}

We have that

\begin{teo}[Quillen]
\label{teo:extension}
The extension $k(A)/k(\eu{g}_A)$ is a normal extension and 
\begin{equation}
 \Aut (k(A)/k(\eu{g}_A))\cong N_G(A)/C_G(A).
\end{equation}
\end{teo}

\begin{proof}
See \cite[Theorem 2.10]{qu:pnil}.
\end{proof}


\subsection{F-isomorphisms and fusion}
\label{ss:cf}
The following lemma is a standard result in commutative algebra.

\begin{lem}
\label{lem:fisomspec}
Let $A$ and $B$ be commutative $\F_p$-algebras and $f\colon A\longrightarrow B$ an $F$-isomorphism. Then $f^*: \Spec (B)\longrightarrow \Spec (A)$ is a homeomorphism. 
\end{lem}

\begin{proof}
Since the kernel of $f$ is 
nilpotent, then for any radical ideal $\eu{a}$ of $A$ one has that $f^{-1} (\sqrt{f(\eu{a})})=\eu{a}$. Since 
for any $x\in B$ there exits $y\in A$ and $n\geq 0$ such that $f(y)=x^{p^n}$, then for any radical ideal $\eu{b}$ of $B$ 
one has that $\sqrt{f (f^{-1}(\eu{b}))}=\eu{b}$. Therefore 
\begin{eqnarray}
\xymatrix{\eu{a}\ar[r]& \sqrt{f(\eu{a})}} \\
\xymatrix{\eu{b}\ar[r]& f^{-1}(\eu{b})}
\end{eqnarray}
is a bijection between the radical ideals of $A$ and the radical ideals of $B$. In particular 
$f^*$ is an isomorphism of varieties.
\end{proof}

We are now ready to prove Theorem A.

\begin{proof}[Proof of Theorem A]
Suppose first that $H$ controls fusion of elementary abelian $p$-subgroups of $G$. Then the embedding functor
\begin{equation}
\label{eq:fisom-1}
\xymatrix{j_{H,G}:\ca{E}_H \ar[r] & \ca{E}_G}
\end{equation}
is an equivalence of categories. Therefore
\begin{equation}
\label{eq:fisom-2}
\xymatrix{H^\bullet (j_{H,G}):H^\bullet(\ca{E}_G) \ar[r] & H^\bullet (\ca{E}_H)}
\end{equation}
is an isomorphism. Consider the commutative diagram  
\begin{equation}
\label{dia:fusF}
\xymatrix{
H^\bullet(G)\ar[d]_{\res^G_H}\ar[r]^{q_G}& H^\bullet(\ca{C}_G)\ar[d]^{H^\bullet(j_{H,G})}\\
H^\bullet(H)\ar[r]^{q_H}& H^\bullet(\ca{C}_H).
}
\end{equation}
By Theorem \ref{thm:strat} and equation \eqref{eq:fisom-2} it follows that $\res^G_H$ is an $F$-isomorphism.

\smallskip

Suppose now that $\res^G_H\colon H^\bullet (G,\F_p)\longrightarrow H^\bullet (H,\F_p)$ is an 
$F$-isomorphism. Then $\res^G_H$ induces an $F$-isomorphism $f\colon H (G)\longrightarrow H (H)$. Consider $A$ and $A^\prime$ two elementary abelian $p$-subgroups of $H$.

\bigskip

\noindent \textit{Subclaim 1:} If $A$ and $A^\prime$ are conjugated in $G$, then they are conjugated in $H$.

\begin{proof}[Subproof]
By Lemma \ref{lem:fisomspec}, $f^*\colon \Spec (H (H))\longrightarrow 
\Spec (H (G))$ provides a bijection between the prime ideals of $H (H)$ and the prime ideals 
of $H (G)$. Furthermore, if $A$ is an elementary abelian $p$-subgroup of 
$H$, then $\eu{g}_A=f^*(\eu{h}_A)$. By Theorem \ref{teo:conjugacion}, if $A$ and $A^\prime$ are conjugated in $G$, then $\eu{g}_A=\eu{g}_{A^\prime}$. In particular $f^* (\eu{h}_A)=\eu{g}_A=\eu{g}_{A^\prime}=f^*(\eu{h}_{A^\prime})$. Therefore 
$\eu{h}_A=\eu{h}_{A^\prime}$. Hence, 
by Theorem \ref{teo:conjugacion}, $A$ and $A^\prime$ are conjugated in $H$.
\end{proof}

\noindent \textit{Subclaim 2:} $N_G(A)=C_G(A)N_H(A)$.

\begin{proof}[Subproof] Since $k(\eu{h}_A)$ is a purely inseparable extension of $k (\eu{g}_A)$, then 
\begin{equation}
 \Aut (k(A)/k(\eu{h}_A))\cong \Aut (k (A)/k (\eu{g}_A)).
\end{equation}
Therefore, by Theorem \ref{teo:extension}, 
$N_H(A)/C_H(A)\cong N_G(A)/C_G(A)$. 
\end{proof}

\noindent \textit{Subclaim 3:} $H$ controls fusion of elementary abelian $p$-subgroups of $G$.

\bigskip

\noindent \textit{Subproof:} Let $A$ be an elementary abelian $p$-subgroup of $H$ and $g\in G$ such that $A^g\leq H$. 
Then, by Subclaim 1 there exists $h\in H$ such that $A^g=A^h$. In particular, by Subclaim 2, 
$gh^{-1}\in N_G(A)=C_G(A).N_H(A)$. Therefore $g\in C_G(A).H$
\end{proof}


\section{A $p$-nilpotency criterion}

\label{s:pnilcriterion}

In this section we will prove our main result Theorem B. To ease the notation we denote by $\Cp$ the class of cyclic groups of order $p$ in case $p$ is odd and cyclic groups of order $2$ and $4$ in case $p=2$. Put $\bp=p$ if $p$ is odd and $\bp=4$ in case $p=2$.  

\begin{teo}
\label{main}
Let $G$ be a finite group and $P$ a Sylow $p$-subgroup of $G$. Then  the following two conditions are equivalent
\begin{enumerate}
\item $G$ is $p$-nilpotent.
\item $P$ controls fusion of $\Cp$-groups. 
\end{enumerate}
\end{teo}

\begin{proof}
It is clear that if $G$ is $p$-nilpotent, then $P$ controls fusion 
of $\Cp$-groups.

Let us show the converse. Using Frobenius 
$p$-nilpotency criterion it is enough to prove that for any subgroup $B$ of $P$ and for any 
$p^\prime$-element $g\in N_G(B)$, then $g$ centralices $B$. The subgroup $B$ 
is contained in  $Z_l(P)$ for some $l\geq 1$ where $Z_l(P)$ denotes the $l$-upper center of $P$. 
We will show by induction on $l$ that $g\in C_G(B)$. 
Suppose first that $B\leq Z(P)$ and consider $a\in B$ such that 
$a^{\bp}=1$. Since $P$ controls 
fusion of $\Cp$-groups, there exists $x\in P$ such that $a^g=a^x$ and since $a\in Z (P)$, then $a^x=a$. 
Hence we have that $g$ centralices all elements of order $p$ ($2$ and $4$ in case $p=2$) in $B$.
Thus, by \cite[Chap. V Lemma 5.12]{hupp1:grp},  $g$ centralices $B$.

For the general case, consider  $B\leq Z_l(P)$ and suppose the assumption to be true for any subgroup contained in $Z_{l-1}(P)$.

\bigskip

\noindent \textit{Subclaim 1:} For $a\in B$ such that $a^{\bp}=1$, we have that $[a,g,g]=1$.

\begin{proof}[Subproof]
We have that  $g$ normalizes the subgroups $K=\langle a\in B\mid a^{\bp}=1\rangle$ and $[K,g]$. We also have that
\begin{equation}
\label{eq-main}
[K,g]=\langle [a,g]^b\mid a,b\in K\ \text{and}\ a^{\bp}=1\rangle .
\end{equation}
Take $a\in B$ such that $a^{\bp}=1$. Since $P$ controls fusion of $\Cp$-groups, there exists $x\in P$  such that $a^g=a^x$. In particular $[a,g]=[a,x]\in Z_{l-1}(P)$. Therefore, 
by \eqref{eq-main}, $[K,g]=Z_{l-1}(P)$. Since $g$ normalices  $[K,g]$ and by induction  hypothesis we have that $[K,g,g]=1$. 
\end{proof}

\noindent \textit{Subclaim 2:} $g\in C_G(B)$.

\bigskip

\noindent \textit{Subproof.}
Take $a\in B$ such that $a^{\bp}=1$ and put $p^e$ the exponent of $B$. Consider the subgroup $H=\langle g,[a,g]\rangle$.  By the Subclaim  1, $\gamma_2(H)=1$. Then, by
 \cite[Chap. III, Theorem 9.4]{hupp1:grp}, we have that $[a,g^{p^e}]= [a,g]^{p^e}=1$.
 Since $g$ is a $p^\prime$-element of $G$, $g$ centralices all elements of order $p$ ($2$ and $4$ in case $p=2$) in $B$. Thus, by
\cite[Chap. V Lemma 5.12]{hupp1:grp},  $g$ centralices $B$. 

This ends the proof.
\end{proof}

As a consequence to this we have the following corollary.

\begin{cor}
\label{corolario1}
Let $G$ a finite group and $P$ a Sylow $p$-subgroup of $G$ such that
\begin{itemize}
\item[1.] $N_G(P)$ controls fusion of $\Cp$-groups and
\item[2.] $N_G(P)=C_G(P).P$.
\end{itemize}
Then $G$ is $p$-nilpotent.  
\end{cor}

\begin{proof}
Let $A$ be a $\Cp$-group  and $g\in G$ such that $A,A^g\leq P$. Since $N_G(P)$ controls fusion  
$\Cp$-groups, one has that $g\in C_G(A).N_G(P)=C_G(A).P$. Then $P$ controls fusion
of $\Cp$-groups  and, by Theorem \ref{main}, $G$ is $p$-nilpotent.
\end{proof}

\section{Some applications}

We now present the first application of Theorem \ref{main}. In \cite{hp:pcen} H-W. Henn 
and S. Priddy proved that if a group $G$ has a Sylow $p$-subgroup $P$ such that 
\begin{itemize}
\item[i)] if $p$ is odd, the elements of order $p$ of $P$ are in the center of $P$ and, if $p=2$, 
the elements of order $2$ and $4$ are in the center of $P$,
\item[ii)] $\Aut (P)$ is a $p$-group,
\end{itemize}
then $G$ is $p$-nilpotent. This implies that "most" finite groups are 
$p$-nilpotent (see \cite{hp:pcen}). The proof of Henn and Priddy is essentially topological. 
In \cite{th} J. Thevenaz gave a group theoretical proof of this result using 
Alperin's Fusion Theorem.  In fact Thevenaz proved that if $G$ satisfies condition i), 
then $N_G(P)$ controls $p$-fusion in $G$. This, together with condition ii) above implies 
that $P$ controls $p$-fusion in $G$ and therefore $G$ is $p$-nilpotent. 
We now give a weaker version of Thevenaz result which
also implies that a group satisfying i) and ii) is $p$-nilpotent.

\begin{pro}
\label{priddy}
Let $G$ be a finite group and $P$ a Sylow $p$-subgroup of $G$. 
Suppose that the elements of order dividing $p$ in $P$ (or $4$ in case $p=2$) 
are in the center of $P$. Then 
$N_G(P)$ controls fusion of $\Cp$-groups.
\end{pro}

\begin{proof}
Let $A$ be a  $\Cp$-group and $g\in G$ such that $A,A^g\leq N_G(P)$. In particular $A^g\leq  P$. Equivalently $A\leq P^{g^{-1}}$. 
Hence, since the elements of $P$ of order $p$ (or $4$ in case $p=2$)  are in the center of $P$, we have that $P,P^{g^{-1}}\leq C_G(A)$. But, since $P$ and $P^{g^{-1}}$ are Syllow $p$-subgroups of $C_G(A)$, there exists $c\in C_G(A)$ such that $P=P^{g^{-1}c}$. Thus $g^{-1}c\in N_G(P)$ and $g\in N_G(P).C_G(A)$.
\end{proof}

\begin{cor}
Let $G$ a finite group and $P$ a Syllow $p$-subgroup of $G$ such that
\begin{itemize}
\item[1.] all elements of order dividing $p$ in $P$ (or $4$ in case $p=2$) 
are in the center of $P$ and
\item[2.] $N_G(P)=P.C_G(P)$.
\end{itemize}
Then $G$ is $p$-nilpotent.
\end{cor}

\begin{proof}
It follows from Proposition \ref{priddy} and Corollary \ref{corolario1}.
\end{proof}

The second application of Theorem \ref{main} is a generalization to $p=2$ of the fact that if the elements of order $p$ of a finite group $G$ are in some upper center of $G$, then $G$ is $p$-nilpotent (see \cite{th:pcenp} and \cite{GSW}).

\begin{cor}
\label{pcentral}
Let $G$ a finite group such that $K=\langle x\in G\mid x^{\bp}=1\rangle\leq Z_n(G)$ 
for some $n\geq 1$ (here $\bp$ means $p$ in case $p$ is odd and $4$ in case $p=2$). 
Then $G$ is $p$-nilpotent.
\end{cor}

\begin{proof}
The subgroup $K$ is nilpotent of class at most $n$, and therefore
a finite $p$-group. Let $p^e$ be the exponent of $K$.
Then, by Hall-Petrescu collection formula (see \cite[Theorem 2.1]{gust:omega}), 
for any $y\in K$ and $x\in G$
\begin{equation}
\label{eq:bas-0}
[y,x^{p^{e+n}}]\in 
\prod_{0\leq i\leq e+n}[K,\overbrace{G,\ldots, G}^{p^i}]^{p^{e+n-i}}=1.
\end{equation}
Therefore one has that $G^{p^{e+n}}\leq C_G(K)$.
Moreover, for any Sylow $p$-subgroup $P$ of $G$ one has
$G=P.G^{p^{e+n}}=P.C_G(K)$.  In particular
$P$ controls fusion of $\Cp$-groups. Hence, by Theorem \ref{main}, $G$ is $p$-nilpotent.
\end{proof}

\bibliography{fusion}
\bibliographystyle{amsplain}
\end{document}